\documentclass[11pt, a4paper]{article}
\usepackage{amsmath,amsfonts,mathrsfs,amssymb,cite,color}
\usepackage[numbers,sort&compress]{natbib}
\usepackage[linkcolor=red,anchorcolor=blue,citecolor=green]{hyperref}

\newtheorem{thm}{\hspace{\parindent}Theorem}[section]%{\hspace{\parindent}Lemma}
\newtheorem{lem}[thm]{\hspace{\parindent}Lemma}
\newtheorem{cor}[thm]{\hspace{\parindent}Corollary}
\newtheorem{rem}[thm]{\hspace{\parindent}Remark}
\newtheorem{prop}[thm]{\hspace{\parindent}Proposition}
\newtheorem{exm}[thm]{\hspace{\parindent}Example}%[section]
\numberwithin{equation}{section}

\newcommand{\C}{\mbox{$\mathbb{C}$}}

\newcommand{\calB}{{\mathcal B}}
\newcommand{\calN}{{\mathcal N}}
\newcommand{\calR}{{\mathcal R}}
\newcommand{\calX}{{\mathcal X}}
\newcommand{\calY}{{\mathcal Y}}

\newcounter{rom}
\newcommand{\xd}[1]{\setcounter{rom}{#1}{\rm (\roman{rom})}}
\title{\Large The invertibility of $2\times2$ operator matrices
\footnotetext{\hspace{-18pt}$^1$Corresponding author.}}

\author{Junjie Huang, \ Junfeng Sun, \ Alatancang Chen, \
Carsten Trunk$^{1}$}

\date{\today}
\begin{document}
\maketitle

\begin{abstract}%establish explicit / smallmatrix
In this paper the properties of right invertible row operators,
i.e., of $1\times 2$ surjective operator matrices are studied.
This investigation is based on a specific space decomposition.
Using this decomposition, we characterize the invertibility of a $2\times2$
 operator matrix. As an application, the invertibility of Hamiltonian
 operator matrices is investigated.
\end{abstract}



\section{Introduction}

The invertibility of a linear operator is one of the most basic problems in operator theory, and, obviously, appears in the study of the linear equation $Tx = y$ with
 a linear operator $T$.

This problem becomes even more involved if one considers the invertibility
of $2\times2$ operator matrices. For this
let $A$, $B$, $C$ and $D$ be bounded linear operators on a Hilbert space. If,
e.g., they are pairwise commutative, then the operator matrix
\begin{equation}\label{MMM}
M=
\begin{pmatrix}
A & B\\
C & D
\end {pmatrix}
\end{equation}
is invertible if and only if  $AD-BC$ is invertible (cf.\ \cite[Problem 70]{Halmos1982}).
If only $C$ and $D$ are commutative, and if,
 in addition, $D$ is invertible, then the operator matrix $M$ is invertible if and only if $AD-BC$ is invertible (cf.\ \cite[Problem 71]{Halmos1982}).
 In fact, the commutativity is essential in the above characterization,
 see \cite[Problem 71]{Halmos1982}. The situation is even more involved
 if $A$ and $D$ are not defined on the same space and, hence, the
 formal expression $AD-BC$ has no meaning.

 In general, there is no complete description of the invertibility of operator matrices in the non-commutative case. But if at least one of the entries $A$ or $D$ of the operator matrix $M$
 is invertible, one can describe the invertibility of $M$ in terms of the Schur complement.
 A similar statement holds also in the case of invertible entries $B$ or $C$.
 Moreover, the Schur complement method can be effectively used also in the
 case where the entries of $M$ are unbounded operators
  under additionally assumptions on the domain of the entries, such as the diagonally (or off-diagonally) dominant or upper (lower) dominant cases,
 see, e.g., the monograph \cite{Tretter2008}.
We also refer to \cite{Kurina2001, Wu-Chen2011} for
 sufficient conditions for
%diagonally (off-diagonally) dominant
nonnegative Hamiltonian operators to have bounded inverses.

However, it is easy to see that there are many invertible $2\times 2$
operator matrices with non invertible entries
$A,B,C$ and $D$ (see, e.g.,
Theorem \ref{special} below).
 Obviously, in such cases, the Schur complement
method is not applicable.

It is the
aim of the present article to give a full characterization for
the invertibility of bounded $2\times 2$ operator matrices. We
do this in the following manner: A necessary condition for the
invertibility of a $2\times 2$
 operator matrix $M$ in~\eqref{MMM} is the fact
that the row operator $(A\ B)$ is right invertible
(that is, the range $\calR((A\ B))$
of the operator $(A\ B)$ covers all of the spaces).
A further necessary condition is  $\calN((A\ B))\ne \{0\}$,
where $\calN((A\ B))$ denotes the kernel of $(A\ B)$ (see
Corollary~\ref{Iliketomoveit} below). This non-zero kernel $\calN((A\ B))$
plays a crucial role. Its projection $P_\calX (\calN((A\ B)))$
onto the first component is a subset of
the kernel of  $P_{\calR(B)^\perp}A$, where $P_{\calR(B)^\perp}$
denotes the orthogonal projection onto $\calR(B)^\perp$. Similarly,
the projection of $\calN((A\ B))$
onto the second component is a subset of $\calN(P_{\calR(A)^\perp}B)$.

Therefore we investigate a right invertible
row operator $(A\ B)$ and choose a decomposition of the space
into six parts which is built out of the subspaces
$\calN(A), \calN(B), \calN(P_{\calR(B)^\perp}A)$
and $\calN(P_{\calR(A)^\perp}B)$.
As a result, we show that the operator $B_2^{-1}
\widetilde{A}_2$ considered as an operator
from $P_\calX (\calN((A\ B)))$ to
$\calN(B)^\perp \ominus \calN(P_{\calR(A)^\perp}B)^\perp$
%$\overline{\calR(A)} \cap \overline{\calR(B)}$
is correctly defined. Here
$\widetilde{A}_2$ ($B_2$) denote the restriction
of $A$ ($B$, respectively) to $\calN(P_{\calR(B)^\perp}A)$
($\calN(B)^\perp \ominus \calN(P_{\calR(A)^\perp}B)^\perp$, respectively).

The main result of the present article is a
full characterization
of the invertibility of a $2\times 2$ matrix operator $M$
in terms of its entries
$A, B, C, D$, or to be more precise,
in terms of the restrictions  $\widetilde{A}_2, B_2,
C_2$ and $D_2$ which are, in some sense, all related to $\calN((A\ B))$:
A $2\times 2$ operator matrix
$M$ is invertible if and only if  the following two statements
are satisfied
\begin{itemize}
\item[\rm (i)] The restriction
$D|_{\calN(B)}$ is left invertible and
\item[\rm (ii)] the operator
$$
C_2 -D_2B_2^{-1}\widetilde{A}_2:
P_\calX (\calN((A\ B))) \to
(\calR(D|_{\calN(B)}))^{\perp}
\mbox{ is one-to-one and surjective.}
$$
\end{itemize}
 Here $C_2$ ($D_2$)  is the restriction of $C$ ($D$, respectively)
to $\calN(P_{\calR(B)^\perp}A)$
($\calN(B)^\perp \ominus \calN(P_{\calR(A)^\perp}B)^\perp$, respectively)
projected onto $(\calR(D|_{\calN(B)}))^{\perp}$.

This characterization is especially helpful if the spaces
$\calN((A\ B))$, $\calN(P_{\calR(B)^\perp}A)$ or
$\calN(P_{\calR(A)^\perp}B)$ are known explicitly,
see, e.g., Theorem~\ref{special} in Section~\ref{Section2}. Moreover, we use
it to derive a characterization for isomorphic row operators
in Section~\ref{Section3}. Finally, in Section~\ref{Section4}
we give an application to Hamiltonian operators.

\section{Main result}\label{Section2}

We always assume that $\calX$ and $\calY$ are complex separable Hilbert spaces.
Let $T$ be a bounded operator between  $\calX$ and $\calY$. We write
$T\in{\mathcal{B}}(\calX,\calY)$ and, if $\calX =\calY$,  $T\in{\mathcal{B}}(\calX)$. The range of $T$ is denoted by
$\calR(T)$, the kernel by $\calN(T)$.
 The term  \textit{isomorphism} is reserved for linear bijections
$T:\calX \to \calY$ that are homeomorphisms, i.e., $T\in {\mathcal{B}}(\calX,\calY)$
and $T^{-1}\in {\mathcal{B}}(\calY,\calX)$.

A subspace in $\calY$ is an operator range if
it coincides with the range of some bounded operator $T\in{\mathcal{B}}(\calX,\calY)$.
The following lemma is from \cite[Theorem 2.4]{Cross1980}.

\begin{lem}\label{lem1}
Let $\calR_{1}$ and $\calR_{2}$ be operator ranges in $\calY$ such that $\calR_{1}+\calR_{2}$ is closed.

\xd{1} If $\calR_{1}\cap{\calR_{2}}$ is closed, then $\calR_{1}$ and $\calR_{2}$ are closed.

\xd{2} If $\calR_{1}$ and $\calR_{2}$ are dense in $\calY$,
then $\calR_{1}\cap{\calR_{2}}$ is dense in $\calY$.
\end{lem}

From \cite[Proposition 2.14, Theorem 2.16]{Brezis2010}, we have the following
basic facts, which are important in the proofs of our main results.

\begin{lem}\label{lem2}
Let $\Omega_1$ and $\Omega_2$ be two closed subspaces in $\calX$.
Then
\[
\Omega_1\cap\Omega_2=(\Omega_1^{\perp}+\Omega_2^{\perp})^{\perp},\quad
\Omega_1^{\perp}\cap{\Omega_2^{\perp}}=(\Omega_1+\Omega_2)^{\perp},
\]
and we further have the following equivalent descriptions:

\xd{1} $\Omega_1+\Omega_2$ is closed;

\xd{2} $\Omega_1^{\perp}+\Omega_2^{\perp}$ is closed;

\xd{3} $\Omega_1+\Omega_2=(\Omega_1^{\perp}\cap{\Omega_2^{\perp}})^{\perp}$;

\xd{4} $(\Omega_1\cap{\Omega_2})^{\perp}=\Omega_1^{\perp}+\Omega_2^{\perp}$.
\end{lem}

As usual, the symbol $\oplus$ denotes the orthogonal sum of
two closed subspaces in a Hilbert space whereas the symbol
 $\dot{+}$ denotes the direct sum of two (not necessarily closed)
 subspaces in a Hilbert space.
 If $\Omega, \Omega_1$ are closed subspaces, $\Omega_1\subset\Omega$,
 we denote by  $\Omega\ominus\Omega_1$ the uniquely determined
 closed subspace  $\Omega_2$ in $\Omega$
  with $\Omega=\Omega_1\oplus\Omega_2$.
%The following lemma is from \cite[Lemma 6]{Markus-Olshevsky1993}.

%\begin{lem}\label{lem4}
%\marginpar{Check}
%Let ${T}\in{{\mathcal{B}}(\calX,\calY)}$. Then for any $\varepsilon>0$,
%there exist orthogonal decompositions $\calX=\calX_{\varepsilon}\oplus{\calX^{\varepsilon}}$
%and $\calY=\calY_{\varepsilon}\oplus{\calY^{\varepsilon}}$ such that
%\[
%\begin{array}{l}
%T(\calX_{\varepsilon})\subset{\calY_{\varepsilon}},\quad \|{Tx}\|_\calY{\leq}\varepsilon{\|x\|_\calX}\
%\mbox{\rm for\ all\ } x\in{\calX_{\varepsilon}},\\
%T(\calX^{\varepsilon})\subset{\calY^{\varepsilon}}\hspace*{-4pt},\quad
%\hspace*{1.5pt}\|{Tx}\|_\calY{\geq}\varepsilon{\|x\|_\calX}\ \mbox{\rm for\ all\ } x\in{\calX^{\varepsilon}}\hspace*{-1pt}.
%\end{array}
%\]
%\end{lem}

The next lemma is well known, see, e.g., \cite[Proposition 1.6.2]{Tretter2008}
or \cite{Nagel1989,Harte1988}.

\begin{lem}\label{lem5}
Let $A\in\calB(\calX), B \in\calB(\calY,\calX), C\in\calB(\calX, \calY)$,
and $D\in\calB(\calY)$. Let $A$ $(B)$ be an isomorphism.
Then the $2\times2$ operator matrix
$$
\begin{pmatrix}
A & B \\C & D
\end{pmatrix}\in\calB(\calX\oplus \calY)
$$
is an isomorphism if and only if $D-CA^{-1}B$
$($resp.\ $C-DB^{-1}A)$
is an isomorphism.
\end{lem}

Recall that an  operator $T\in{\mathcal{B}}(\calX,\calY)$ is called
right invertible if there exists an operator
$S\in{\mathcal{B}}(\calY,\calX)$ with $TS=I_\calY$, where $I_\calY$
stands for the identity mapping in $\calY$. Hence, if $T$ is right
invertible then it is surjective.
Conversely, if $T\in{\mathcal{B}}(\calX,\calY)$
then the restriction
$T|_{\mathcal{N}(T)^\perp}$ maps $\mathcal{N}(T)^\perp$ onto $\calR(T)$ and, if
$\calR(T)=\calY$, then $T|_{\mathcal{N}(T)^\perp}: \mathcal{N}(T)^\perp \to \calY$
is an isomorphism. Then with
\begin{equation}\label{Erfurt1}
S:=\left( \begin{matrix} 0  \\ \left(T|_{\mathcal{N}(T)^\perp}\right)^{-1}  \end{matrix} \right)
:\calY\to \mathcal{N}(T)\oplus \mathcal{N}(T)^\perp
\end{equation}
considered as an operator in ${\mathcal{B}}(\calY,\calX)$ we see that
$T$ is right invertible. This shows the equivalence of (i)-(iii) in
the following (well-known) lemma.
\begin{lem}\label{Herzz}
For $T\in{\mathcal{B}}(\calX,\calY)$ the following assertions are equivalent.
\begin{itemize}
\item[\rm (i)] The operator $T$ is right invertible.
\item[\rm (ii)] $\calR(T)= \calY$.
\item[\rm (iii)] The operator $T|_{\mathcal{N}(T)^\perp}$
considered as an operator from
$\mathcal{N}(T)^\perp$ into
$\calY$
is an isomorphism.
\item[\rm (iv)] There exists an isomorphism $U\in {\mathcal{B}}(\calY)$ such that
$UT$ is a right invertible operator.
\end{itemize}
\end{lem}

\indent{\textit{Proof.~~}}
It remains to show the equivalence of (iv) with (i)-(iii). Choose $U=I_\calY$ and we see that
(i) implies (iv). Conversely, let $U\in {\mathcal{B}}(\calY)$  be an isomorphism.
If $UT$ is right invertible, then by (ii) $\calR(UT)=\calY$. As $\calR(T)=\calR(UT)$, again
(ii) shows that $T$ is right invertible.
\hfill $\Box$\\

Similarly, $T\in{\mathcal{B}}(\calX,\calY)$ is called
left invertible if there exists an operator
$S\in{\mathcal{B}}(\calY,\calX)$ with $ST=I_\calX$.
 Hence, if $T$ is left
invertible then it is injective
and for a sequence $(y_n)$ in $\calR(T)$ with $y_n\to y$ as $n\to \infty$
we find $(x_n)$ with $Tx_n=y_n$ and
$$
x_n = STx_n = Sy_n \to Sy \quad \mbox{and} \quad
y_n =Tx_n \to TSy,
$$
which shows the closedness of $\calR(T)$.

Conversely, if $\calN(T)=\{0\}$ and $\calR(T)$ is closed,
then $T$ considered as an operator from $\calX$ into $\calR(T)$
is  an isomorphism and its inverse $T^{-1}$ acts from $\calR(T)$ into
$\calX$.  Then with
\begin{equation}\label{Erfurt2}
S:=\left( 0 \ \ T^{-1}  \right): \calR(T)^\perp\oplus \calR(T) \to \calX,
\end{equation}
considered as an operator in ${\mathcal{B}}(\calY,\calX)$, we see that
$T$ is left invertible. We collect these statements in the following lemma,
where the equivalence of (i)-(iii) follows from the above considerations
and the equivalence of (i)-(iii) with (iv) is obvious.

\begin{lem}\label{Lefin}
For $T\in{\mathcal{B}}(\calX,\calY)$ the following assertions are equivalent.
\begin{itemize}
\item[\rm (i)] The operator $T$ is left invertible.
\item[\rm (ii)] $\calN(T)= \{0\}$ and $\calR(T)$ is closed.
\item[\rm (iii)] The operator $T$ considered as an operator from
$\calX$ into $\calR(T)$ is an isomorphism.
\item[\rm (iv)] There exists an isomorphism $V\in {\mathcal{B}}(\calX)$ such that
$TV$ is a left invertible operator.
\end{itemize}
\end{lem}

\begin{rem}\label{Erfurt3}
The following observation for $T\in{\mathcal{B}}(\calX,\calY)$
follows immediately from the Lemmas \ref{Herzz} and \ref{Lefin}.
If $T$ is right invertible, then there exists a left invertible operator
$S\in  {\mathcal{B}}(\calY,\calX)$ (cf.\ \eqref{Erfurt1})
with $TS=I_\calY$ and $\calR(S)=\mathcal{N}(T)^\perp$. If
$T$ is left invertible,  then there exists a right invertible operator
$S\in  {\mathcal{B}}(\calY,\calX)$ (cf.\ \eqref{Erfurt2})
with $ST=I_\calX$.
\end{rem}

For the orthogonal
projection onto a closed subspace $\Omega$ in some Hilbert space we
shortly write $P_\Omega$.

\begin{thm}\label{thm1}
Let ${A}\in{\mathcal{B}}(\calX)$ and ${B}\in{\mathcal{B}}(\calY,\calX)$
 and assume that the row operator $(A\ B) \in{\mathcal{B}}(\calX \oplus\calY,\calX)$
is right invertible.  Then $\calX$ admits the decomposition
\begin{equation}\label{katholisch0}
\calX =(\calR(A)^{\bot}\dot{+}{\calR(B)^{\perp}})
\oplus{\overline{\calR(A)}\cap{\overline{\calR(B)}}}
\end{equation}
and the space $\calX \oplus \calY$ admits the decomposition
\begin{equation}\label{katholisch1}
\calX \oplus \calY =
\calX_1\oplus\calX_2\oplus\calX_3\oplus\calY_3\oplus\calY_2\oplus\calY_1,
\end{equation}
where
\begin{equation}\label{katholisch2}
\begin{array}{l}
\calX_1:={{{\mathcal{N}}}(A)},\quad
\calX_2:={{\mathcal{N}}}(A)^{\perp}\ominus
{{{\mathcal{N}}(P_{\calR(B)^{\perp}}A)^{\perp}}},\quad
\calX_3:={{\mathcal{N}}(P_{\calR(B)^{\bot}}A)^{\bot}};\\
\calY_1:={{{\mathcal{N}}}(B)},\quad
\calY_2:={{\mathcal{N}}}(B)^{\perp}\ominus{{{\mathcal{N}}(P_{\calR(A)^{\perp}}B)^{\perp}}},\quad
\calY_3:={{\mathcal{N}}(P_{\calR(A)^{\perp}}B)^{\perp}}.
\end{array}
\end{equation}
The row operator $(A \  B)$ from  $\calX \oplus \calY$ into
$\cal X$ admits the following  representation with respect to the
decompositions \eqref{katholisch0} and \eqref{katholisch1}
\begin{equation}\label{katholisch3}
\begin{pmatrix}
0 & 0 & 0 & B_{3}& 0 & 0 \\
0 & 0 & A_{3} & 0 & 0 & 0 \\
0 & A_{2} & A_0 & B_0 & B_{2} & 0
\end{pmatrix},
\end{equation}
where
\begin{equation*}
\begin{array}{l}
A_0\in{\mathcal{B}}\left(\calX_3,\overline{\calR(A)}\cap\overline{\calR(B)}\right),
\; A_2\in{\mathcal{B}}\left(\calX_2,\overline{\calR(A)}\cap\overline{\calR(B)}\right),
\; A_3\in{\mathcal{B}}\left(\calX_3,\calR(B)^\perp\right);\\[1ex]
B_0\in{\mathcal{B}}\left(\calY_3,\overline{\calR(A)}\cap\overline{\calR(B)}\right),
\; B_2\in{\mathcal{B}}\left(\calY_2,\overline{\calR(A)}\cap\overline{\calR(B)}\right),
\; B_3\in{\mathcal{B}}\left(\calY_3,\calR(A)^\perp\right).
\end{array}
\end{equation*}
Then the operators $A_3$ and $B_3$ are isomorphisms
and the row operator $(A_{2}~B_{2}): \calX_2\oplus\calY_2\rightarrow {\overline{\calR(A)}\cap{\overline{\calR(B)}}}$
is right invertible and
\begin{equation}\label{katholisch4}
\overline{\calR(A_2)}=\overline{\calR(A)}
\cap\overline{\calR(B)}=\overline{\calR(B_2)}.
\end{equation}
\end{thm}

\indent{\textit{Proof.~~}}
\textit{Step 1. We prove \eqref{katholisch0}--\eqref{katholisch3}.}

The row operator  $( A \ B): \calX\oplus \calY\to \calX$ is right invertible
and we have with Lemma~\ref{Herzz}
\begin{equation}\label{alpha.alpha}
\calR(A)+\calR(B)=\calX.
\end{equation}
We claim
\begin{equation}\label{Eins}
P_{\calR(A)^{\perp}}(\calR(B))={\calR(A)^{\perp}}.
\end{equation}
To see this, it suffices to show the inclusion $P_{\calR(A)^{\perp}}(\calR(B))\supset{\calR(A)^{\perp}}$.
Let $x\in\calR(A)^{\perp}$. Then there exist $x_1\in\calR(A)$ and $x_2\in\calR(B)$ such that
$x=x_1 + x_2$, so $x=P_{\calR(A)^{\perp}}x_2\in P_{\calR(A)^{\perp}}(\calR(B))$. This proves
the claim. Similarly, we obtain
\begin{equation}\label{Zwei}
P_{\calR(B)^{\perp}}(\calR(A))={\calR(B)^{\perp}}.
\end{equation}
\iffalse
Using Lemma \ref{lem3}, we deduce
\[
\begin{array}{l}
AA^{\ast}|_{\calR(B)^{\bot}}=AA^{\ast}|_{{{\mathcal{N}}}(B^{\ast})}+BB^{\ast}|_{{{\mathcal{N}}}(B^{\ast})}
\geq\delta {I|_{{{\mathcal{N}}}(B^{\ast})}}=\delta{I|_{\calR(B)^{\bot}}},
\\
BB^{\ast}|_{{{\mathcal{R}}}(A)^{\bot}}=AA^{\ast}|_{{{\mathcal{N}}}(A^{\ast})}+BB^{\ast}|_{{{\mathcal{N}}}(A^{\ast})}\geq\delta
{I|_{{{\mathcal{N}}}(A^{\ast})}}=\delta{I|_{\calR(A)^{\bot}}}.
\end{array}
\]
Then we immediately obtain the inclusions $\calR(A)\supseteq{\calR(B)^{\bot}}$ and $\calR(B)\supseteq{\calR(A)^{\bot}}$.
\fi
Moreover, by \eqref{alpha.alpha}, we have
\[
\{0\}=\calX^\perp=(\overline{\calR(A)}+\overline{\calR(B)})^{\perp}=
\overline{\calR(A)}^{\perp}\cap{\overline{\calR(B)}^{\perp}}
\]
and also the sum $\overline{\calR(A)}+\overline{\calR(B)}$ is closed.
By Lemma \ref{lem2} (iv) it follows that
\[
\left(\overline{\calR(A)}\cap{\overline{\calR(B)}}\right)^{\perp}
=\overline{\calR(A)}^{\perp}+\overline{\calR(B)}^{\perp}.
\]
To sum up, we have the space decomposition \eqref{katholisch0}.
As $\mathcal N(A)\subset
\mathcal N (P_{\calR(B)^{\perp}}A)$, we have
$\mathcal N (P_{\calR(B)^{\perp}}A)^\perp \subset\mathcal N(A)^\perp$.
Analogously we see
$\mathcal N (P_{\calR(A)^{\perp}}B)^\perp \subset\mathcal N(B)^\perp$
and, hence,  decomposition \eqref{katholisch1} follows.

For $x\in \calX_3^\perp =\mathcal{N}(P_{\calR(B)^{\bot}}A)$ we have
\begin{equation*}
Ax= \left(I-P_{\calR(B)^{\bot}}\right)Ax= P_{\overline{\calR(B)}}Ax.
\end{equation*}
Hence, $x\in \mathcal{N}(P_{\calR(B)^{\bot}}A)$ if and only if
\begin{equation}\label{Stern}
Ax\in \overline{\calR(B)}.
\end{equation}
Similarly, $y\in \mathcal{N}(P_{\calR(A)^{\bot}}B)$ if and only if
$By\in \overline{\calR(A)}$.
Therefore, if $x_2\in \calX_2$ ($y_2\in \calY_2$), then it follows that $x_2\in \mathcal{N}(P_{\calR(B)^{\bot}}A)$ (resp.\ $y_2\in \mathcal{N}(P_{\calR(A)^{\bot}}B)$) and, by \eqref{Stern}
\begin{equation}\label{Herz}
Ax_2\in \overline{\calR(B)} \quad
(\mbox{resp.\ } By_2\in \overline{\calR(A)}).
\end{equation}
Then the zero entries in \eqref{katholisch3}
 follow from the fact that $Ax=0$ for $x\in
\mathcal N (A)$, $By=0$ for $y\in
\mathcal N (B)$, $Ax\in \calR(A)$,  $By\in \calR(B)$, and \eqref{Herz}.

\ \\
\indent
\textit{Step 2. We show that $(A_{2}~B_{2})$ is right invertible.}

We have $\mathcal N (A) \subset \mathcal{N}(P_{\calR(B)^{\bot}}A)$,
$\mathcal N (B) \subset \mathcal{N}(P_{\calR(A)^{\bot}}B)$
and by \eqref{alpha.alpha} and \eqref{katholisch0} we see that
$A_3$ and $B_3$ are isomorphisms.
Thus, there exists an isomorphism
$U\in\calB(
(\calR(A)^{\bot}\dot{+}{\calR(B)^{\perp}})
\oplus{\overline{\calR(A)}\cap{\overline{\calR(B)}}})
$
$$
U:= \begin{pmatrix}
1 & 0 & 0  \\
0 & 1 & 0\\
-B_0B_{3}^{-1} & -A_0A_3^{-1}  & 1
\end{pmatrix}
$$
such that
\begin{equation*}\label{equ21}
U\begin{pmatrix}
0 & 0 & 0 & B_{3}& 0 & 0 \\
0 & 0 & A_{3} & 0 & 0 & 0\\
0 & A_{2} & A_0 & B_0 & B_{2} & 0
\end{pmatrix}=\begin{pmatrix}
0 & 0 & 0 & B_{3}& 0 & 0 \\
0 & 0 & A_{3} & 0 & 0 & 0 \\
0 & A_{2} & 0 & 0 & B_{2} & 0
\end{pmatrix}.
\end{equation*}
As $(A \ B)$ is right invertible, Lemma \ref{Herzz} shows that $(A_{2}~B_{2}): \calX_2\oplus\calY_2\rightarrow {\overline{\calR(A)}\cap{\overline{\calR(B)}}}$
is right invertible.

\ \\
\indent
\textit{Step 3. We  show \eqref{katholisch4}.}

By definition, we have $\calR(A_2)\subset \overline{\calR(A)}
\cap\overline{\calR(B)}$ and $\calR(B_2)\subset \overline{\calR(A)}
\cap\overline{\calR(B)}$. We will only show $\overline{\calR(A)}
\cap\overline{\calR(B)}\subset \overline{\calR(B_2)}$.
The proof for  $\overline{\calR(A)}
\cap\overline{\calR(B)}\subset \overline{\calR(A_2)}$ is the same and,
hence, we omit this proof.

Let $z\in \overline{\calR(A)}
\cap\overline{\calR(B)}$. Then there exists a sequence $(z_n)$ in
$\calR(B)$ which converges to $z$. By the  block representation
 \eqref{katholisch3} for $B$ we find $z_{1,n}$ in $\calR(A)^\perp$
 and $z_{3,n} \in \overline{\calR(A)}
\cap\overline{\calR(B)}$ with
\begin{equation} \label{Advent2a}
z_n= z_{1,n} +z_{3,n}, \quad n\in \mathbb N,
\end{equation}
where we have
\begin{equation} \label{Advent2}
z_{1,n} = B_3y_{3,n} \quad \mbox{and} \quad z_{3,n}= B_0 y_{3,n} + B_2y_{2,n}
\quad \mbox{for }
n\in \mathbb N
\end{equation}
 for some  $y_{2,n}\in \calY_2$ and $y_{3,n}\in  \calY_3$.
The convergence of $(z_n)$ implies the convergence of $(z_{1,n})$ to some
$z_1\in\calR(A)^\perp$ and of $(z_{3,n})$ to some $z_3 \in \overline{\calR(A)}
\cap\overline{\calR(B)}$,
$$
z=z_1 + z_3.
$$
The vectors $z$ and $z_3$ belong to $ \overline{\calR(A)}$, thus $z_1
\in  \overline{\calR(A)}$ and $z_1=0$ follows. Therefore $(B_3y_{3,n} )$ in
\eqref{Advent2} converges to zero. The fact that $B_3$ is an
isomorphism implies $y_{3,n} \to 0$ as $n\to\infty$. We conclude
$$
z=z_3= \lim_{n\to\infty} z_{3,n} = \lim_{n\to\infty} B_2y_{2,n}
$$
and $z\in  \overline{\calR(B_2)}$ follows. Relation \eqref{katholisch4} is proved.
\hfill $\Box$\\

The following proposition will be used in the proof of the main result.

\begin{prop}\label{thm11}
Let ${A}\in{\mathcal{B}}(\calX)$ and ${B}\in{\mathcal{B}}(\calY,\calX)$
and let the row operator $(A\ B) \in{\mathcal{B}}(\calX \oplus\calY,\calX)$
be right invertible. The following assertions are equivalent.
\begin{itemize}
\item[\rm (i)]
 $\calR(B)$ is closed.
 \item[\rm (ii)]
  $P_{\calX}({{\mathcal{N}}}((A \ B)))$
is a closed subspace in $\calX$.
\item[\rm (iii)]
$\calR(B_2)$ is closed.
\end{itemize}
\end{prop}

\indent{\textit{Proof.~~}}
Let $\calR(B)$ be closed. We have
$$
P_{\calX}({{\mathcal{N}}}((A  \ B)))
=\{x\in \calX: Ax\in{\mathcal{R}(A)\cap{\mathcal{R}(B)}}\}=
\{x\in \calX: Ax\in\mathcal{R}(B)\}
$$
and $P_{\calX}({{\mathcal{N}}}((A  \ B)))$
is the pre-image of $\calR(B)$ under $A$, and, hence,  it
is a closed subspace and (ii) holds.

If $P_{\calX}({{\mathcal{N}}}((A \ B)))$ is closed, then also
 $$
 \Omega:=P_{\calX}({{\mathcal{N}}}((A \ B)))\cap \mathcal N(A)^\perp=
\{x\in \calX : x\in{{{\mathcal{N}}}(A)^{\bot}}, \; Ax\in{\calR(A)\cap{\calR(B)}}\}
 $$
is closed. Decompose $x\in\Omega$ with respect to the decomposition,
cf.\ Theorem \ref{thm1},
$\calX= \calX_1 \oplus \calX_2 \oplus \calX_3$ as $x=x_1+x_2+x_3$ with
$x_j\in\calX_j$ for $j=1,2,3$. Then $x_1=0$ and for some $y\in \calY$ we have
$Ax=By$. Decompose $y$ with respect to
$\calY= \calY_1 \oplus \calY_2 \oplus \calY_3$
(cf.\ Theorem \ref{thm1}) as $y=y_1+y_2+y_3$ with
$y_j\in\calY_j$ for $j=1,2,3$. Relation \eqref{katholisch3} shows
$$
Ax= A
\left( \begin{matrix} 0 \\ x_2\\ x_3 \end{matrix} \right)=
\left( \begin{matrix} 0 \\ A_3x_3\\ A_2x_2+A_0x_3 \end{matrix} \right)=
\left( \begin{matrix} B_3y_3 \\ 0\\ B_0y_3+B_2y_2 \end{matrix} \right)=
B
\left( \begin{matrix} y_3  \\ y_2\\ y_1 \end{matrix} \right) =By
$$
and, as $A_3$ is an isomorphism, we obtain $x_3=0$. Therefore
$\Omega\subset \calX_2$ and we write
$$
\calX_2= \Omega \oplus (\calX_2\ominus{\Omega}).
$$
By Theorem \ref{thm1} $(A_{2}~B_{2})$ is right invertible and we obtain with Lemma~\ref{Herzz}
\[
A_{2}(\calX_2\ominus{\Omega})+B_{2}(\calY_2)={\overline{\calR(A)}\cap{\overline{\calR(B)}}},\quad
A_{2}(\calX_2\ominus{\Omega})\cap{B_{2}(\calY_2)}=\{0\}.
\]
Thus, using Lemma \ref{lem1}, we deduce that $A_{2}(\calX_2\ominus{\Omega})$ and
$\calR(B_{2})$ are closed.

Assume that (iii) holds. Then,
by \eqref{katholisch4}, the operator $B_2$ is an isomorphism.
Let $z\in\overline{\calR(B)}$. Then there exists a sequence $(z_n)$ in
$\calR(B)$ which converges to $z$. By the  block representation
 \eqref{katholisch3} for $B$ we find $z_{1,n}$ in $\calR(A)^\perp$
 and $z_{3,n} \in \overline{\calR(A)}
\cap\overline{\calR(B)}$ such that
 \eqref{Advent2a} and \eqref{Advent2} hold
 for some  $y_{2,n}\in \calY_2$ and $y_{3,n}\in  \calY_3$.
The convergence of $(z_n)$ implies the convergence of $(z_{1,n})$ to some
$z_1\in\calR(A)^\perp$ and of $(z_{3,n})$ to some $z_3 \in \overline{\calR(A)}
\cap\overline{\calR(B)}$, $z=z_1 + z_3$. As the operators
$B_2$ and $B_3$ (cf.\ Theorem \ref{thm1}) are isomorphisms, we have
$$
y_{3,n} \to B_3^{-1} z_1 \quad y_{2,n} \to -B_2^{-1}B_0 B_3^{-1} z_1 +
B_2^{-1}z_3 \quad \mbox{as } n\to \infty.
$$
Thus, with  \eqref{katholisch3},
$$
B
\left( \begin{matrix}  B_3^{-1} z_1 \\
-B_2^{-1}B_0 B_3^{-1} z_1 +B_2^{-1}z_3 \\ 0 \end{matrix} \right)=
\left( \begin{matrix} z_1 \\ 0\\ z_3 \end{matrix} \right)
=z,
$$
and $z\in\calR(B)$.
\hfill $\Box$\\

\begin{lem}\label{lll}
Let ${A}\in{{\mathcal{B}}(\calX)}$, ${B}\in{{\mathcal{B}}(\calY,\calX)}$ and
 assume that the row operator $(A\ B) \in{\mathcal{B}}(\calX \oplus\calY,\calX)$
is right invertible. Let $A_2$ and $B_2$ be
  as in Theorem \ref{thm1}.
  Then $B_2$ considered as an operator from $\mathcal Y_2$
  to $\mathcal R(B_2)$ is one-to-one and has an inverse
  $B_2^{-1} :\mathcal R(B_2) \to \mathcal Y_2$. Define
  $$
   \widetilde{A}_2:=(0\ A_2):
   \mathcal X_1\oplus \mathcal X_2
   \to{\overline{\calR(A)}\cap{\overline{R(B)}}}.
$$
Then  $\widetilde{A}_2|_{P_{\calX}({{\mathcal{N}}}((
A \ B)))}$ maps to $\mathcal R(B_2)$ and the operator
$$
B_2^{-1}\widetilde{A}_2|_{P_{\calX}({{\mathcal{N}}}((
A \ B)))}: P_{\calX}({{\mathcal{N}}}((
A \ B))) \to \mathcal Y_2
$$
 is correctly defined.

 If $\calR(B)$ is closed, then
$B_{2}$ is an isomorphism and we have
$$
\mathcal X_1\oplus \mathcal X_2=
{{\mathcal{N}}(P_{\calR(B)^{\bot}}A)}=P_{\calX}({{\mathcal{N}}}( (
A \ B )))
$$
and the operator
\begin{equation}\label{llla}
B_2^{-1}\widetilde{A}_2: {{\mathcal{N}}(P_{\calR(B)^{\bot}}A)} \to \mathcal Y_2
\end{equation}
 is correctly defined.
\end{lem}

\indent{\textit{Proof.~~}}
As $\mathcal Y_2\subset \mathcal{N}(B)^\perp$ the operator
$B_2$ is one-to-one, hence its inverse
$B_2^{-1} :\mathcal R(B_2) \to \mathcal Y_2$ exists. From
\begin{equation}\label{chef}
P_{\calX}({{\mathcal{N}}}((A \ B)))=
\{x\in \calX : Ax\in{\calR(A)\cap{\calR(B)}}\}
\subset \{x\in \calX : Ax\in\overline{\calR(B)} \}
\end{equation}
we conclude
$$
P_{\calX}({{\mathcal{N}}}((A \ B)))\subset{{\mathcal{N}}(P_{\calR(B)^{\bot}}A)}
= \mathcal X_1\oplus \mathcal X_2.
$$
Moreover, we decompose $x\in P_{\calX}({{\mathcal{N}}}((A \ B)))$
 with respect to the decomposition
$\calX= \calX_1 \oplus \calX_2 \oplus \calX_3$
(cf.\ Theorem \ref{thm1}) as $x=x_1+x_2+x_3$ with
$x_j\in\calX_j$ for $j=1,2,3$. Then $x_3=0$ and for some $y\in \calY$ we have
$Ax=By$. Decompose $y$ with respect to
$\calY= \calY_1 \oplus \calY_2 \oplus \calY_3$
(cf.\ Theorem \ref{thm1}) as $y=y_1+y_2+y_3$ with
$y_j\in\calY_j$ for $j=1,2,3$. Relation \eqref{katholisch3} shows
$$
Ax= A
\left( \begin{matrix} x_1 \\ x_2\\ 0\end{matrix} \right)=
\left( \begin{matrix} 0 \\ 0\\ A_2x_2 \end{matrix} \right)=
\left( \begin{matrix} B_3y_3 \\ 0\\ B_0y_3+B_2y_2 \end{matrix} \right)=
B
\left( \begin{matrix} y_3  \\ y_2\\ y_1 \end{matrix} \right) =By
$$
and, as $B_3$ is an isomorphism, we obtain $y_3=0$ and $A_2x_2=B_2y_2$.
Thus $\widetilde{A}_2x\in \mathcal R(B_2)$ for $x\in
P_{\calX}({{\mathcal{N}}}((A \ B)))$ and $B_2^{-1}\widetilde{A}_2|_{P_{\calX}({{\mathcal{N}}}((
A \ B)))}$ is correctly defined.
If $\calR(B)$ is closed, then by
Proposition \ref{thm11} also $\calR(B_2)$ is closed and
by \eqref{katholisch4} we see that $B_2$ is an isomorphism.
Moreover, from \eqref{chef} we see in this case
$ \mathcal X_1\oplus \mathcal X_2=
{{\mathcal{N}}(P_{\calR(B)^{\bot}}A)}=P_{\calX}({{\mathcal{N}}}( (
A \ B )))$ and \eqref{llla} follows.
\hfill $\Box$\\

The following theorem is the main result. It provides a full
characterization of isomorphic $2\times 2$ operator matrices
in terms of their entries.

\begin{thm}\label{thm2}
Let ${A}\in{{\mathcal{B}}(\calX)}$, ${B}\in{{\mathcal{B}}(\calY,\calX)}$.
 Assume that the row operator $(A\ B) \in{\mathcal{B}}(\calX \oplus\calY,\calX)$
is right invertible and, hence,  adopt the notions $A_2$, $B_2$,
and $\mathcal X_j$, $\mathcal Y_j$, $j=1,2,3$,
as in Theorem \ref{thm1} and $\widetilde{A}_2$ as in Lemma \ref{lll}. Let
${C}\in{{\mathcal{B}}(\calX,\calY)}$ and ${D}\in{{\mathcal{B}}(\calY)}$.
Define the operator matrix $M$ by
\[
M= \begin{pmatrix}
A & B \\
C & D
\end{pmatrix}.
\]
Define the operator $B_2^{-1}\widetilde{A}_2|_{P_{\calX}({{\mathcal{N}}}((
A \ B)))}$  as in Lemma \ref{lll} and define
$$
C_2:=P_{(\calR(D|_{\calN(B)}))^{\perp}} C|_{\mathcal X_1\oplus \mathcal X_2}:
X_1\oplus \mathcal X_2 \to (\calR(D|_{\calN(B)}))^{\perp}
$$
and
$$
D_2:=P_{(\calR(D|_{\calN(B)}))^{\perp}} D|_{\calY_2}:
\calY_2 \to (\calR(D|_{\calN(B)}))^{\perp}.
$$
Then $M$ is an isomorphism if and only if
 the following two statements are satisfied:
\begin{itemize}
\item[\rm (i)]
  The restriction $D|_{\mathcal{N}(B)}:\mathcal{N}(B)\to\calY$ is left invertible.
\item[\rm (ii)]
The operator
$$
\left.  \left(C_{2}-D_{2}B_{2}^{-1}
\widetilde{A}_{2}\right)\right|_{P_{\calX}({{\mathcal{N}}}((
A \ B)))}:P_{\calX}({{\mathcal{N}}}((A \ B)))
\to (\calR(D|_{\calN(B)}))^{\perp}
$$
 is one-to-one and surjective.
\end{itemize}
\end{thm}

\indent{\textit{Proof.~~}}
Let $M$ be an isomorphism.
Then the row operator
 $(A \ B): \calX\times\calY\to\calX$ is right invertible,
see Lemma \ref{Herzz}, and
the column operator
$\left(\begin{smallmatrix}B \\ D\end{smallmatrix}\right):
\calY \to \calX\times\calY$ is injective.
Moreover, if the range of $\left(\begin{smallmatrix}B \\ D\end{smallmatrix}\right)$
is not closed then there exists a sequence $(y_n)$ in $\calY$
with $\|y_n\|=1$, $n\in \mathbb N$, and
$\left(\begin{smallmatrix}B \\ D\end{smallmatrix}\right)y_n\to 0$
 as $n\to \infty$. But this implies
 $M\left(\begin{smallmatrix}0 \\ y_n\end{smallmatrix}\right)\to 0$,
 a contradiction as $M$ is assumed to be an isomorphism. Therefore
 the column operator
$\left(\begin{smallmatrix}B \\ D\end{smallmatrix}\right)$
is left invertible, cf.\ Lemma \ref{Lefin}.

Now let $z\in \overline{\calR(D|_{\calN(B)})}$. Then, there exists $z_n\in \calN(B)$ such that
$Dz_n\to z$ as $n\to\infty$, and we further have
\[
\begin{pmatrix}B \\ D\end{pmatrix}z_n
=\begin{pmatrix}0 \\ Dz_n\end{pmatrix}
\to\begin{pmatrix}0 \\ z\end{pmatrix},
\]
which together with Lemma \ref{Lefin} implies
\[
\begin{pmatrix}B \\ D\end{pmatrix}x
=\begin{pmatrix}0 \\ z\end{pmatrix}
\]
for some $x\in\calN(B)$, and hence $D|_{\calN(B)}x=z$.
This proves that $\calR(D|_{\calN(B)})$ is closed, hence,
$D|_{{\mathcal{N}}(B)}$ is left invertible by Lemma \ref{Lefin}
and (i) is proved.

As $\calR(D|_{\calN(B)})$ is a closed subspace in $\calY$,
we decompose $\calY$,
\begin{equation}\label{zerleg12}
\calY=(\calR(D|_{\calN(B)}))^{\perp}\oplus
\calR(D|_{\calN(B)}).
\end{equation}

Similar to the proof of Theorem \ref{thm1}, $M$ %=\left(\begin{smallmatrix} A & C \\ D & B \end{smallmatrix}\right)$
as an operator from ${{\mathcal{N}}(P_{\calR(B)^{\bot}}A)}\oplus\calX_3
\oplus\calY_3\oplus\calY_2\oplus\calY_1$ into
$$
(\calR(A)^{\bot}\dot{+}{\calR(B)^{\perp}})
\oplus{\overline{\calR(A)}\cap{\overline{R(B)}}}
\oplus(\calR(D|_{\calN(B)}))^{\perp}\oplus
\calR(D|_{\calN(B)})
$$
has the following block representation
\begin{equation}\label{eq+}
M=\left(
\begin{matrix}
0 & 0 & B_{3}& 0 & 0 \\
0 & A_{3} & 0 & 0 & 0 \\
\widetilde{A}_{2} & A_0 & B_0 & B_{2} & 0 \\
C_{2} & C_{3} & D_{1} & D_{2} & 0\\
C_{4} & C_{5} & D_{3} & D_{4} & D_{5}
 \end{matrix}
 \right).
\end{equation}
By Theorem \ref{thm1},
 $A_{3}$ and $B_{3}$  are isomorphisms. Additionally, as $M$ is an isomorphism,
 $D_5$ is also an isomorphism.
Then there exist isomorphisms
\begin{equation*}
\begin{array}{l}
U\in\calB\left((\calR(A)^{\bot}\dot{+}{\calR(B)^{\perp}})\oplus{\overline {\calR(A)}\cap{\overline{R(B)}}}\oplus(\calR(D|_{\calN(B)}))^{\perp}
\oplus\calR(D|_{\calN(B)})\right),\\[1ex]
V\in\calB\left(\mathcal{N}(P_{\calR(B)^\bot} A)\oplus\calX_3\oplus\calY_3\oplus\calY_2\oplus\calY_1\right)
\end{array}
\end{equation*}
with
$$
U:= \begin{pmatrix}
1 & 0 & 0 & 0 & 0 \\
0 & 1 & 0 & 0 & 0\\
-B_0B_{3}^{-1} & -A_0A_3^{-1}  & 1 & 0 & 0\\
-D_1B_{3}^{-1} & -C_3A_3^{-1}  & 0 & 1 & 0\\
0 & 0  & 0 & 0 & 1
\end{pmatrix},
$$
$$
V:= \begin{pmatrix}
1 & 0 & 0 & 0 & 0 \\
0 & 1 & 0 & 0 & 0\\
0 & 0 & 1 & 0 & 0\\
0 & 0 & 0 & 1 & 0\\
-D_5^{-1}C_{4} & -D_5^{-1}C_{5} & -D_5^{-1}D_{3} & -D_5^{-1}D_{4} & 1
\end{pmatrix}
$$
such that
\begin{equation}\label{eq++}
UMV=\left(
\begin{matrix}
0 & 0 & B_{3}& 0 & 0 \\
0 & A_{3} & 0 & 0 & 0 \\
\widetilde{A}_{2} & 0 & 0 & B_{2} & 0 \\
C_{2} & 0 & 0 & D_{2} & 0\\
0 & 0 & 0 & 0 & D_{5}
 \end{matrix}
 \right).
\end{equation}
Thus, $M$ is an isomorphism if and only if
\begin{equation}\label{eq+++}
\Delta:=\left(\begin{matrix}\widetilde{A}_{2} & B_{2} \\C_{2} & D_{2}\end{matrix}\right):
{{\mathcal{N}}(P_{\calR(B)^{\bot}}A)}\oplus{{\calY_2}}\rightarrow({\overline{\calR(A)}
\cap{\overline{\calR(B)}}})\oplus(\calR(D|_{\calN(B)}))^{\perp}
\end{equation}
is an isomorphism.

{\it Case 1:} $\calR(B)$ is closed. In this case, from Lemma \ref{lll},
$B_{2}:\calY_2\rightarrow{\overline{\calR(A)}\cap{\overline{\calR(B)}}}$ is an isomorphism and
$B_{2}^{-1}\widetilde{A}_{2}:
{{\mathcal{N}}(P_{\calR(B)^{\bot}}A)} \to \mathcal Y_2$ is correctly defined,
see Lemma \ref{lll}.
According to Lemma \ref{lem5}, $\Delta$ is an isomorphism if and only if
%$\left(\begin{smallmatrix}0 & B_{2} \\C_{2}-D_{2}B^{-1}_{2}A_{2} &  0\end{smallmatrix}\right)$ is invertible, and if and only if
\[
C_{2}-D_{2}B_{2}^{-1}\widetilde{A}_{2}:
{{\mathcal{N}}(P_{\calR(B)^{\bot}}A)}\rightarrow (\calR(D|_{\calN(B)}))^{\perp}
\]
is an isomorphism. By Lemma \ref{lll} ${{\mathcal{N}}(P_{\calR(B)^{\bot}}A)}=
P_{\calX}({{\mathcal{N}}}((A \ B)))$ and (ii) is satisfied.

{\it Case 2:} $\calR(B)$ is not closed. By Proposition \ref{thm11} also
$\calR(B_2)$ is not closed which implies $\dim \calR(B_2)=\infty$ and
$\dim \calY_2 = \infty$. The dimension does not change when we close
a subspace, therefore we conclude from \eqref{katholisch4}
\begin{equation}\label{NullNull}
\dim \overline{\calR(A)}
\cap\overline{\calR(B)}=\dim \overline{\calR(B_2)}=
\dim \calR(B_2)=\infty.
\end{equation}
By Theorem \ref{thm1} $(A_{2} \  B_{2})$ is right invertible,
\eqref{katholisch4} and Lemma \ref{lem1} imply
$$
\overline{\calR(A_2)\cap \calR(B_2)} = \overline{\calR(A)}
\cap \overline{\calR(B)}.
$$
Obviously, $\calR(A_2)\cap \calR(B_2) \subset \calR(A)
\cap \calR(B)$ and we obtain $\overline{\calR(A)}
\cap \overline{\calR(B)} \subset \overline{\calR(A)
\cap \calR(B)}$. Thus
\begin{equation*}\label{Pfeil}
\overline{\calR(A)}
\cap \overline{\calR(B)} = \overline{\calR(A)\cap \calR(B)}.
\end{equation*}
From this and from $\calR(A)\cap \calR(B) \subset
\calR(A)\cap \overline{\calR(B)} \subset
\overline{\calR(A)}\cap \overline{\calR(B)}$ we conclude
with \eqref{NullNull}
\begin{equation}\label{HerzHerz}
\infty =\dim \overline{\calR(A)\cap \calR(B)} =
\dim \calR(A)\cap \calR(B) = \dim \calR(A)\cap \overline{\calR(B)}.
\end{equation}
We will use \eqref{HerzHerz} to show
\begin{equation}\label{KaroKaro}
\dim{{{\mathcal{N}}}((\widetilde{A}_{2} \  B_{2}))}
=\dim{{{\mathcal{N}}(P_{\calR(B)^{\bot}}A)}}.
\end{equation}
For this we consider
\begin{equation}\label{Punkte}
{{\mathcal{N}}}((A \ B))=\left\{
\left(\begin{smallmatrix}x  \\ 0  \end{smallmatrix}\right) :
x\in {{\mathcal{N}}}(A)\right\}
\oplus{\left\{\left(\begin{smallmatrix}y  \\ z  \end{smallmatrix}\right):
y\in{{{\mathcal{N}}}(A)^{\bot}},  Ay=-Bz\right\}}
\end{equation}
and
\begin{equation*}
\mathcal{N}(P_{\calR(B)^{\bot}}A)= \mathcal{N}(A) \oplus
\left\{x : x \in  \mathcal{N}(A)^\perp, Ax \in \overline{\calR(B)}\right\}.
\end{equation*}
As $A$ restricted to $\calN(A)^\perp$ is injective, we obtain with
\eqref{HerzHerz}
\begin{equation*}
 \begin{split}
\dim {\left\{\left(\begin{smallmatrix}y  \\ z  \end{smallmatrix}\right):
y\in{{{\mathcal{N}}}(A)^{\bot}},  Ay=-Bz\right\}}& =
\dim  \calR(A)\cap \calR(B)= \dim  \calR(A)\cap \overline{\calR(B)}\\
&= \dim \left\{x : x \in  \mathcal{N}(A)^\perp, Ax \in \overline{\calR(B)}\right\}.
\end{split}
\end{equation*}
Therefore
\begin{equation*}
\dim {{\mathcal{N}}}((A \ B)) = \dim \mathcal{N}(P_{\calR(B)^{\bot}}A)
\end{equation*}
and with \eqref{eq++} we obtain $\dim{{{\mathcal{N}}}((\widetilde{A}_{2} \  B_{2}))}
=\dim{{{\mathcal{N}}(P_{\calR(B)^{\bot}}A)}}$, hence \eqref{KaroKaro} is
proved. Two separable Hilbert spaces of the same dimension are
unitarily equivalent, therefore there exists a left invertible operator
\begin{equation}\label{beta.beta}
\left(\begin{matrix}G \\H\end{matrix}\right):\calY_2\rightarrow{{\mathcal{N}}(P_{\calR(B)^{\bot}}A)}
\oplus\calY_2 \mbox{ with range } {{\mathcal{N}}}((\widetilde{A}_{2} \ B_{2})).
\end{equation}
Since $\calX_1\oplus \calX_2
 = \mathcal{N}(P_{\calR(B)^{\bot}}A)$
 and by Theorem \ref{thm1} and Lemma \ref{lll}
 $(\widetilde{A}_{2} \ B_{2}):
 \mathcal{N}(P_{\calR(B)^{\bot}}A)\oplus \calY_2\to
 \overline{\calR(A)}\cap \overline{\calR(B)}$ is a right invertible operator. Then, see Remark \ref{Erfurt3}, there exists a left invertible operator
\begin{equation}\label{beta.beta3}
\left(\begin{matrix}E\\F \end{matrix}\right):
\overline{\calR(A)}\cap \overline{\calR(B)}
\rightarrow{{{\mathcal{N}}
(P_{\calR(B)^{\bot}}A)}\oplus{{\calY_2}}}
\end{equation}
such that
\begin{equation}\label{beta.beta2}
\widetilde{A}_{2}E+B_{2}F=I_{{\overline{\calR(A)}\cap{\overline{\calR(B)}}}}
 \quad \mbox{with } \calR\left(\left(\begin{matrix}E\\F \end{matrix}\right)\right)=
  ({{\mathcal{N}}}((\widetilde{A}_{2} \ B_{2})))^\perp
\end{equation}
Define
\begin{equation}\label{beta.beta4}
W=\left(\begin{matrix}E & G \\F & H\end{matrix}\right):
\overline{\calR(A)}\cap \overline{\calR(B)}
\oplus\calY_2
\rightarrow\mathcal{N}(P_{\calR(B)^{\bot}}A)\oplus\calY_2.
\end{equation}
As $\left(\begin{smallmatrix}G   \\ H  \end{smallmatrix}\right)$
and $\left(\begin{smallmatrix}E  \\ F  \end{smallmatrix}\right)$
are left invertible and from \eqref{beta.beta} and \eqref{beta.beta2}
we obtain easily that $W$ is an isomorphism. We have
\begin{equation}\label{allday0}
\Delta W=\left(\begin{matrix}I_{{\overline{\calR(A)}\cap{\overline{\calR(B)}}}} & 0 \\C_{2}E+D_{2}F & C_{2}G+D_{2}H\end{matrix}\right).
\end{equation}
As $M$
is an isomorphism, $\Delta$ is an isomorphism (see  \eqref{eq+++})
and the operator $C_{2}G+D_{2}H: \calY_2\to (\calR(D|_{\calN(B)}))^{\perp}$
is an isomorphism. Moreover, the operator $B_2$ considered as an operator
from $\calY_2$ to $\calR(B_2)$ is one-to-one and has an inverse,
see Lemma \ref{lll}. From $\widetilde{A}_{2}G+B_{2}H=0$ we conclude
$-B_{2}^{-1} \widetilde{A}_{2}G=H$ and
\begin{equation}\label{allday1}
C_{2}G+D_{2}H=(C_{2}-D_{2}B_{2}^{-1} \widetilde{A}_{2})G.
\end{equation}
Therefore, $C_{2}-D_{2}B_{2}^{-1} \widetilde{A}_{2}: \calR (G)
\to (\calR(D|_{\calN(B)}))^{\perp}$ is one-to-one with range equal
to $(\calR(D|_{\calN(B)}))^{\perp}$. From
\begin{equation}\label{allday2}
\begin{split}
\calR(\left(\begin{smallmatrix}G   \\ H  \end{smallmatrix}\right))
&={{\mathcal{N}}}( (\widetilde{A}_{2} \ B_{2} ))\\
&=\left(\begin{smallmatrix}{{\mathcal{N}}}(A)  \\ 0  \end{smallmatrix}\right)
\oplus{\left\{\left(\begin{smallmatrix}x   \\ y  \end{smallmatrix}\right):
x\in{{{\mathcal{N}}}(A)^{\bot}},  y\in{{{\mathcal{N}}}(B)^{\bot}},Ax=-By\right\}}\\
&=\mathcal{N}((A \ B)),
\end{split}
\end{equation}
see \eqref{Punkte},
it follows that $\calR(G)=P_{\mathcal{X}}(\mathcal{N}( (A \ B )))$
and (ii) is shown.
\ \\

Now let us assume that (i) and (ii) hold.
Then $\calR(D|_{\calN(B)})$ is a closed subspace and  $\calY$
admits a decomposition as in \eqref{zerleg12} and we obtain
 the representation of $M$ as in \eqref{eq+}, where $A_3$, $B_3$ and $D_5$
 are isomorphisms. Then, taking the same $U$ and $V$ as above, we obtain
  the relation \eqref{eq++}. Moreover, if $\Delta$ in (\ref{eq+++}) is an isomorphism, then $M$ is an isomorphism.

 If $\calR(B)$ is closed, then from Lemma \ref{lll},
$B_{2}:\calY_2\rightarrow{\overline{\calR(A)}\cap{\overline{\calR(B)}}}$ is an isomorphism and
$B_{2}^{-1}\widetilde{A}_{2}:
{{\mathcal{N}}(P_{\calR(B)^{\bot}}A)} \to \mathcal Y_2$ is correctly defined.
Moreover, Lemma \ref{lll}, ${{\mathcal{N}}(P_{\calR(B)^{\bot}}A)}=
P_{\calX}({{\mathcal{N}}}((A \ B)))$.
Then, by (ii),
\[
C_{2}-D_{2}B_{2}^{-1}\widetilde{A}_{2}:
{{\mathcal{N}}(P_{\calR(B)^{\bot}}A)}\rightarrow (\calR(D|_{\calN(B)}))^{\perp}
\]
is an isomorphism and  according to Lemma \ref{lem5}, $\Delta$ is
an isomorphism and, hence, $M$ is an isomorphism.

If $\calR(B)$ is not closed, then as above, we define the operators
$G$, $H$, $E$, $F$, and $W$ as in \eqref{beta.beta},  \eqref{beta.beta3}, \eqref{beta.beta2}, and \eqref{beta.beta4}. Moreover,
the operator $W$ in \eqref{beta.beta4} is an isomorphism and also
 \eqref{allday1} and  \eqref{allday2} hold.
By \eqref{allday2} $\calR(G)=P_{\mathcal{X}}(\mathcal{N}( (A \ B )))$ and
as $B_2$ is one-to-one,
we see that the operator $G$ in \eqref{beta.beta} is one-to-one.
Hence, together with (ii), the operator
$(C_{2}-D_{2}B_{2}^{-1} \widetilde{A}_{2})G
: \calY_2\to (\calR(D|_{\calN(B)}))^{\perp}$
is one-to-one with range equal to $(\calR(D|_{\calN(B)}))^{\perp}$.
Therefore, by  \eqref{allday1}, $C_{2}G+D_{2}H$ is an
isomorphism and, by \eqref{allday0} and as $W$  is an isomorphism,
also $\Delta$  is an isomorphism. Therefore, see
(\ref{eq+++}), $M$ is an isomorphism.
\hfill $\Box$\\

Finally, we consider the following special case.
\begin{thm}\label{special}
 Let $A,B,C,D\in \mathcal B(\calX)$ and let $\calX', \calX''$
 be closed subspaces of $\calX$ with
 $$
 \calX = \calX' \oplus  \calX''
 $$
 such that
 $$
 \calR(A)=  \calX', \quad \calN(A) =  \calX'', \quad
 \calR(B)= \calX'', \quad \mbox{and} \quad \calN(B)=  \calX'.
 $$
 Moreover assume that the restriction $D|_{\calX'}:
 \calX' \to \calX$ is left invertible.
 Then  the $2\times 2$ operator matrix $M$,
\[
M= \begin{pmatrix}
A & B \\
C & D
\end{pmatrix},
\]
 is an isomorphism if and only if
$$ C_2:= P_{(\calR(D|_{\calX'}))^{\perp}}
\left. C\right|_{\calX''}: \calX''
\to (\calR(D|_{\calX'}))^{\perp}
$$
 is an isomorphism.

 In particular, if, in addition, $\calR(B) \ne \{0\}$
 and the operator $D|_{\calX'}:
 \calX' \to \calX$ is an isomorphism, then for every operator
 $C\in {{\mathcal{B}}(\calX)}$  the $2\times 2$ operator matrix $M$
 is not an isomorphism.
\end{thm}

{\textit{Proof.}}
Denote by $P_\calX$ the orthogonal projection in $\calX \oplus \calX$
onto the first component. Then
$$
P_\calX (\calN((A\ B))) = \calN(A) = \calX''.
$$
Moreover, we have ${\mathcal{N}}(P_{\calR(B)^{\perp}}A)^{\perp}
= {\mathcal{N}}(P_{\calX'}A)^{\perp} = \calN(A)^\perp$ and
$\overline{\calR(A)} \cap \overline{\calR(B)} = \calX'
\cap \calX'' = \{0\}$. Then the space $\calX_2$
in Theorem \ref{thm1} equals zero and the operators $A_2$
and $\widetilde{A}_{2}$ in Theorem \ref{thm2} are zero.
Then the statements of Theorem \ref{special} follow from
 Theorem \ref{thm2}.
\hfill $\Box$\\

\section{A characterization of isomorphic row operators}
\label{Section3}

In this section let $A,B,C,D$ and $M$ be as in Theorem~\ref{thm2}.
In the following we use Theorems~\ref{thm1} and \ref{thm2}
to characterize the case of an isomorphic row operator $(A\ B)$
and to derive a necessary condition for
$M$ to be an isomorphism.

\begin{prop}\label{Kuerbis}
Let ${A}\in{\mathcal{B}}(\calX)$ and ${B}\in{\mathcal{B}}(\calY,\calX)$.
The row operator $(A\ B) \in{\mathcal{B}}(\calX \oplus\calY,\calX)$
is an isomorphism $($i.e.\ $(A\ B)$ is left and right invertible$)$ if
and only if the following two statements are satisfied:
\begin{itemize}
\item[\rm (i)] $\calN(A) = \calN(B) = \{0\}.$
\item[\rm (ii)] $\calR(A) = \calR(B)^\perp$, $\calR(B)=\calR(A)^\perp$.
\end{itemize}
\end{prop}

{\textit{Proof.}}
If (i) and (ii) hold, then $Ax+By=0$ for some $x\in \calX$,
$y\in \calY$ implies $Ax=-By\in\calR(B)$. By (ii), $Ax=0$
and, hence, $By=0$ follows. Then (i) implies $x=y=0$ and
$\calN((A\ B))=\{0\}$. Moreover, we have with (ii)
$$
\calR ((A\ B)) \subset \calR(A) + \calR(B)=
\calR(A) + \calR(A)^\perp =\calX
$$
and the row operator $(A\ B)$ is an isomorphism.

For the contrary let the row operator $(A\ B)$ be an isomorphism.
If for some $x\in \calX$ we have $Ax=0$ then
$(A\ B) \left(\begin{smallmatrix}x \\ 0\end{smallmatrix}\right)=0$
and, as $\calN(A\ B) =\{0\}$,
 $x=0$ follows. That is, $\calN(A) = \{0\}$ and,
similarly, we see $\calN(B) = \{0\}$. This shows (i). In order to
show (ii) let $x\in \overline{\calR(A)}\cap \overline{\calR(B)}$ and
assume $x\ne 0$. Then there exists sequences $(x_n)$ in $\calX$
and $(y_n)$ in $\cal Y$ such that $(Ax_n)$ and $(B y_n)$ converge both to $x$
with $\liminf_{n\to\infty} \|x_n\| >0$ and
$\liminf_{n\to\infty} \|y_n\| >0$.
But then $(A\ B)  \left(\begin{smallmatrix}x_n \\ -y_n\end{smallmatrix}\right)= Ax_n-By_n$ tends to zero
and $\calR((A\ B))$ is not closed, a contradiction.
This shows
\begin{equation}\label{HTag}
\overline{\calR(A)}\cap \overline{\calR(B)} = \{0\}.
\end{equation}
As $x\in \calN(P_{\calR(B)^\perp}A)$ if and only if
$Ax \in \overline{\calR(B)}$ (see also \eqref{Stern}),
we conclude with $\calN(A)=\{0\}$ and \eqref{HTag}
$$
\calN(P_{\calR(B)^\perp}A)=\{0\}.
$$
In the same way we obtain from \eqref{HTag} and $\calN(B)=\{0\}$
that $\calN(P_{\calR(A)^\perp}B)=\{0\}$. Then  for
the spaces $\calX_1, \calX_2, \calX_3, \calY_1, \calY_2, \calY_3$ from
Theorem~\ref{thm1} we conclude
$$
\calX_1=\{0\}, \quad \calX_2=\{0\}, \quad \calX_3=\calX, \quad
\calY_1=\{0\}, \quad \calY_2=\{0\}, \quad
\mbox{and} \quad \calY_3=\calY
$$
and the row operator $(A\ B)$
admits  a representation according to Theorem~\ref{thm1} with
respect to the decompositions
$\calX \oplus \calY$
and $\calX =\calR(A)^{\bot}\dot{+}{\calR(B)^{\perp}}$
of the form
\begin{equation*}
\begin{pmatrix}
 0 & B_{3} \\
 A_{3} & 0 \\
\end{pmatrix},
\end{equation*}
where
$A_3\in{\mathcal{B}}\left(\calX,\calR(B)^\perp\right)$ and
$B_3\in{\mathcal{B}}\left(\calY,\calR(A)^\perp\right)$ are isomorphisms.
This shows (ii).
\hfill $\Box$\\

\begin{exm}
Let $\calX=\calY=\ell^2(\mathbb N)$ and
consider the following operators $A$ and $B$ in $X$:
\begin{equation*}
 A(x_n)_{n\in\mathbb N}:=(x_1,0,x_2,0\dots)\quad \mbox{and} \quad
 B(x_n)_{n\in\mathbb N}:=(0,x_1,0,x_2\dots).
\end{equation*}
Then the  row operator $(A\ B)$ satisfies {\rm (i)} and {\rm (ii)} of
Proposition~\ref{Kuerbis}
and, hence, $(A\ B)$ is an isomorphism.
\end{exm}

As a consequence, we derive the following condition
for $M$ to be an isomorphism.
\begin{cor}\label{Iliketomoveit}
Let ${A}\in{{\mathcal{B}}(\calX)}$, ${B}\in{{\mathcal{B}}(\calY,\calX)}$,
${C}\in{{\mathcal{B}}(\calX,\calY)}$ and ${D}\in{{\mathcal{B}}(\calY)}$.
If
$$
\calY \ne \{0\} \quad \mbox{and} \quad \calN((A\ B))=\{0\}
$$
then the operator matrix $M$
\[
M= \begin{pmatrix}
A & B \\
C & D
\end{pmatrix}
\]
is not a isomorphism.
\end{cor}

{\textit{Proof.}}
If $M$ is an isomorphism, then
as noted in the proof of Theorem~\ref{thm2}, the row
operator $(A\ B)$ is right invertible. Assume
$\calN((A\ B))=\{0\}$. Then $(A\ B)$ is an isomorphism,
and, by Proposition~\ref{Kuerbis}, $\calN(B)
=\{0\}$. Hence, we obtain $(\calR(D|_{\calN(B)}))^{\perp}=\calY$
and (ii) in Theorem~\ref{thm2}
cannot be true unless $\calY=\{0\}$. Therefore,
either $\calY=\{0\}$ or $\calN((A\ B))\ne\{0\}$ holds.
\hfill $\Box$\\

\section{Application to  Hamiltonian operators}\label{Section4}

In this section we consider the special case
of Hamiltonian operators, i.e.,
in the situation of Theorem \ref{thm2},
$\calX=\calY$, the operators
$B,C$ are self-adjoint and $D=-A^*$.
Under these assumptions, Theorem~\ref{thm2}
takes the following simple form.

\begin{thm}\label{cor1} Let $A,B,C\in{{\mathcal{B}}(\calX)}$.
Assume that the row operator $(A\ B) \in{\mathcal{B}}(\calX \oplus\calX,\calX)$
is right invertible  and   that ${B}$ and ${C}$
 are self-adjoint operators in $\calX$, i.e.\ $B=B^*$ and $C=C^*$.
Adopt the notions $A_2$, $B_2$,
and $\mathcal X_j$, $\mathcal Y_j$, $j=1,2,3$,
as in Theorem \ref{thm1} and $\widetilde{A}_2$ as in Lemma \ref{lll}.
Define the operator $B_2^{-1}\widetilde{A}_2|_{P_{\calX}({{\mathcal{N}}}((
A \ B)))}$  as in Lemma \ref{lll} and define
$$
C_2:=P_{\calN(P_{\calR(B)^\perp}A)} C|_{\mathcal X_1\oplus \mathcal X_2}:
X_1\oplus \mathcal X_2 \to \calN(P_{\calR(B)^\perp}A)
$$
and
$$
(-A^*)_2:=-P_{\calN(P_{\calR(B)^\perp}A)} A^*|_{\calY_2}:
\calY_2 \to \calN(P_{\calR(B)^\perp}A).
$$
Then the Hamiltonian operator
\[
H=\begin{pmatrix}
A & B \\ C & -A^*
\end{pmatrix}
\]%\in{{\mathcal{B}(X\oplus{X})}}$
is an isomorphism if and only if
\begin{itemize}
\item[\rm (i)]
the operator
  $$
  \left.\left(C_{2}-(-A^*)_2B_{2}^{-1}
  \widetilde{A}_{2}\right)\right|_{P_{\calX}({{\mathcal{N}}}((A \ B)))}:
P_{\calX}({{\mathcal{N}}}((A \ B)))
\to \calN(P_{\calR(B)^\perp}A)
$$
 is one-to-one and surjective.
\end{itemize}

If in this case we have, in addition, that
$\calR(B)$ is closed, then $C_{2}-(-A^*)_2B_{2}^{-1}
  \widetilde{A}_{2}\in\calB(\calN(P_{\calR(B)^\perp}A))$
%:P_{\calX}({\mathcal{N}}( ( A \ B ) ))\rightarrow P_{\calX}({\mathcal{N}}( ( A \ B ) ))$
is an isomorphism.
\end{thm}

{\textit{Proof.}}
By assumption, the row operator $(A\ B)$
is right invertible, hence (see Lemma \ref{Herzz}) its range is closed and
$\calR(A)+\calR(B)={\calX}$. The same applies to
$(B\ -A)$ and thus its adjoint,
$$
(B\ -A)^* = \left(\begin{matrix}B \\-A^*\end{matrix}\right),
$$
has a closed range and is one-to-one.
Let $z\in \overline{\calR(-A^*|_{\calN(B)})}$.
Then, there exists $z_n\in \calN(B)$ such that
$-A^*z_n\to z$ as $n\to\infty$, and we further have
\[
\begin{pmatrix}B \\ -A^*\end{pmatrix}z_n
=\begin{pmatrix}0 \\ -A^*z_n\end{pmatrix}
\to\begin{pmatrix}0 \\ z\end{pmatrix},
\]
which together with the closedness of the range of $(B\ -A)^*$ implies
\[
\begin{pmatrix}B \\ -A^*\end{pmatrix}x
=\begin{pmatrix}0 \\ z\end{pmatrix}
\]
for some $x\in\calN(B)$, and hence $-A^*|_{\calN(B)}x=z$.
This proves that $\calR(-A^*|_{\calN(B)})$ is closed
and (i) in Theorem \ref{thm2} is satisfied for $D=-A^*$.

Next, we verify
\begin{equation}\label{Rennsteigkreuzung}
(\calR(-A^*|_{\calN(B)}))^{\perp}
=\calN(P_{\calR(B)^\perp}A).
\end{equation}
Indeed, if $x\in(\calR(-A^*|_{\calN(B)}))^{\perp}$, we have $(-Ax,y)=(x,-A^{\ast}y)=0$ for every
$y\in{{\mathcal{N}}(B)}$, hence $-Ax \in \calN(B)^\perp$,
which together with the self-adjointness of $B$ deduces $Ax\in\overline{\calR(B)}$, and hence $x\in\calN(P_{\calR(B)^\perp}A)$; while if $x\in\calN(P_{\calR(B)^\perp}A)$, then $Ax\in{\overline{\calR(B)}}$, and hence we  have for $y\in{{\mathcal{N}}(B)}$ that
$(x,-A^{\ast}y)=(-Ax,y)=0$,
i.e., $x\in(\calR(-A^*|_{\calN(B)}))^{\perp}$.

Now the equivalence of (i) and the fact that $H$ is
an isomorphism follows from \eqref{Rennsteigkreuzung} and Theorem \ref{thm2}. The
additional statement in the case of a closed range
of $B$ follows from Lemma \ref{lll}.
\hfill $\Box$\\

\section*{Acknowledgements}

Junjie Huang gratefully acknowledges the support by the National Natural
Science Foundation of China (No.\ 11461049), and the Natural Science Foundation of Inner Mongolia (No.\ 2013JQ01).

\subsection*{Contact information}

{\bf Junjie Huang}

School of Mathematical Sciences, Inner Mongolia University

010021 Hohhot, P.R.\ China

huangjunjie@imu.edu.cn

\noindent
{\bf Junfeng Sun}

School of Mathematical Sciences, Inner Mongolia University

010021 Hohhot, P.R.\ China

sunjunfeng20099@163.com

\noindent
{\bf Alatancang Chen}

School of Mathematical Sciences, Inner Mongolia University

010021 Hohhot, P.R.\ China

alatanca@imu.edu.cn

\noindent
{\bf Carsten Trunk }

Institut f\"ur  Mathematik,  Technische Universit\"{a}t Ilmenau

Postfach 100565, D-98684 Ilmenau,  Germany

carsten.trunk@tu-ilmenau.de


\begin{thebibliography}{10}
\bibitem{Brezis2010}{H.\ Brezis}, Functional Analysis, Sobolev Spaces and Partial Differential Equations,
Springer, New York, 2010.

%\bibitem{Cao-Guo-Meng2006}{X.\ Cao, M.\ Guo, B.\ Meng},
%Semi-Fredholm spectrum and Weyl`s theorem for operator matrices, Acta Math.\ Sin.\ (Engl.\ Ser.), 22 (2006), 169--178.

%\bibitem{Conway1990}{J.B.\ Conway}, A Course in Functional Analysis, Second Edition, Springer-Verlag, New York, 1990.

\bibitem{Cross1980}{R.W.\ Cross}, On the continuous linear image of a Banach space, J.\ Austral.\ Math.\ Soc.\
(Ser.\ A) 29 (1980), 219--234.


%\bibitem{Douglas1966}{R.G.\ Douglas}, On majorization, factorization, and range inclusion of operators on Hilbert space, Pro.\ Amer.\ Math.\ Soc., 17 (1966), 413--415.

%\bibitem{Du-Pan1994}{H.\ Du, J.\ Pan}, Perturbation of spectrums of $2\times2$ operator matrices, Proc.\ Amer.\ Math.\ Soc.\ 121 (1994), 761--766.

%\bibitem{Gohberg-Goldberg1981}{I.\ Gohberg, S.\ Goldberg}, Basic operator theory, Birkh\"{a}user, Boston, 1981.

%\bibitem{Gohberg-Goldberg-Kaashoek1993}I.\ Gohberg, S.\ Goldberg, M.A.\ Kaashoek, Classes of Linear Operators II, Birkh${\rm\ddot{a}}$user-Verlag, Basel, 1993.


%\bibitem{Hai-Chen2010} G.\ Hai, A.\ Chen, Completion problems and spectra for operator matrices, 517 (2010), 31--78.

\bibitem{Halmos1982}{P.R.\ Halmos}, A Hilbert Space Problem Book, Second Edition, Springer, New York, 1982.

%\bibitem{Han-Lee-Lee2000}{J.K.\ Han, H.Y.\ Lee, W.Y.\ Lee}, Invertible completions of $2\times2$ upper trianglar operator matrices, Proc.\ Amer.\ Math.\ Soc., 128 (2000), 119--123.

\bibitem{Harte1988}{R.\ Harte},
Invertibility and singularity of operator matrices,
Proc.\ R.\ Ir.\ Acad.\ 88A (1988), 103--118.


%\bibitem{Hwang-Lee2001}{I.S.\ Hwang, W.Y.\ Lee}, The boundedness below of $2\times2$ upper trianglar operator matrices, Integr.\ Equ.\ Oper.\ Th.\ 39 (2001), 267--276.

\bibitem{Kurina2001}{G.A.\ Kurina},
Invertibility of nonnegatively Hamiltonian operators
in a Hilbert space, Differential Equations 37 (2001), 880--882.

%\bibitem{LangerMarkusMatsaevTretter2001}H.\ Langer, A.\ Markus, V.\ Matsaev, C.\ Tretter,  A new concept for block operator matrices: the quadratic numerical range, Linear Algebra Appl.\ 330 (2001), 89--112.

%\bibitem{Li2004}Y.\ Li, The Perturbations on Spectra of Operator Matrix, Master Degree Thesis, Shanxi Normal University, 2004.

%\bibitem{Markus-Olshevsky1993}{A.S.\ Markus, V.R.\ Olshevsky},
%Complete controllability and spectrum assignment in infinite dimensional spaces,
%Integral Equations Operator Theory 17 (1993), 109--121.

\bibitem{Nagel1989}{R.\ Nagel},
Towards a matrix theory for unbounded operator matrices,
Math.\ Z.\ 201 (1989), 57--68.


%\bibitem{Ren-Du2004}F.\ Ren, H.\ Du, On the invertible completions of operator partial matrices, Journal of Central China Normal University (Natural Science Edition), 38 (2004), 15--16.

%\bibitem{Ren-Huang2006}F.\ Ren, J.\ Huang, On the invertible completions of operator partial matrices, Journal of Northwest University (Natural Science Edition), 36 (2006), 645--648.

%\bibitem{Takahashi1995}K.\ Takahashi, Invertible completions of operator matrices, Integr.\ Equ.\ Oper.\ Th.\ 21 (1995), 355--361.\


\bibitem{Tretter2008}
C.\ Tretter,
Spectral Theory of  Block Operator Matrices
and Applications, Imperial College Press, London,  2008.

%\bibitem{Tretter2009} C.\ Tretter, Spectral inclusion for unbounded block operator matrices, J.\ Funct.\ Anal.\ 256 (2009), 3806--3829.

%\bibitem{TretterWagenhofer2003} C.\ Tretter, M.\ Wagenhofer, The block numerical range of an $n\times n$ block operator matrix, SIAM.\ J.\ Matrix Anal.\ \& Appl.\ 24 (2003), 1003¨C1017.

%\bibitem{Taylor-Lay1980}{A.E.\ Taylor, D.C.\ Lay}, Introduction to Functional Analysis, 2nd ed., Joln Wiley \& Sons, New York, 1980.

%\bibitem{Weidmann1980}{J.\ Weidmann}, Linear operator in Hilbert spaces, Springer-Verlag, New York, 1980.

\bibitem{Wu-Chen2011}
D.\ Wu, A.\ Chen, Invertibility of nonnegative Hamiltonian operator with
unbounded entries, J.\ Math.\ Anal.\ Appl.\ 373 (2011), 410--413.

\end{thebibliography}
\end{document}